\documentclass[12pt]{amsart}

\author{Yutaka Hemmi}
\address{Department of Mathematics and Information Science,
  Faculty of Science, Kochi University, Kochi 780--8520}
\email{hemmi@math.kochi-u.ac.jp}
\title{Retractions of H-spaces}
\newtheorem{thm}{Theorem}[section]
\newtheorem{lemma}[thm]{Lemma}
\numberwithin{equation}{section}
\newcommand{\bmath}[1]{\mbox{\boldmath$#1$}}

\begin{document}
\begin{abstract}
Stasheff showed that if a map between $H$-spaces is an $H$-map,
then the suspension of the map is extendable to a map between
projective planes of the $H$-spaces.
Stahseff also proved the converse under the assumption
that the multiplication of the target space of the map is homotopy associative.
We show by giving an example that the assumption of homotopy associativity
of the multiplication of the target space is necessary to show the converse.
We also show an analogous fact for maps between $A_n$-spaces.
\end{abstract}
\maketitle

\section{Introduction}

Let $X$ and $Y$ be $H$-spaces, and $f \colon X \to Y$ a map.
Stasheff \cite{Stas61a} showed that if $f$ is an $H$-map,
then it's suspension $\Sigma f \colon \Sigma X \to \Sigma Y$ is extendable to
a map $P_2f \colon P_2X \to P_2Y$ between projective planes $P_2X$ and $P_2Y$ of $X$ and $Y$, respectively.
He also showed the converse under the assumption
that the multiplication $\mu_Y$ of $Y$ is homotopy associative.
It has not been known if the converse holds without the assumption of the homotopy associativity of
$\mu_Y$.
In this paper 
we show by giving an example that the assumption of homotopy associativity
of $\mu_Y$ is necessary to show the converse.

Our example is the retraction $r \colon J(X) \to X$ for an $H$-space $X$.
Here, $J(X)$ is the reduced power space of $X$, which is defined as an identification space
of $\bigcup_{i\ge 1} X^i$.
Then the map $r$ is defined by
\[
r([x_1,\dots,x_i])=(\cdots((x_1\cdot x_2)\cdot x_3)\cdots )\cdot x_i,
\]
where $[x_1,\dots,x_i]$ is the class of $(x_1,\dots,x_i) \in X^i$ and 
$x\cdot y$ denotes the multiplication of $x$ and $y$.
Our result is stated as follows.

\begin{thm}
\label{thm:h}
For any $H$-space $X$, there is an extension $P_2r \colon P_2J(X) \to P_2X$ of
$\Sigma r \colon \Sigma J(X) \to \Sigma X$.
\end{thm}

Stasheff showed the following
\begin{thm}[\cite{Stas61a}]
\label{thmm}
The retraction $r$ is an $H$-map if and only if the multiplication
of $X$ is homotopy associative.
\end{thm}

Thus in particular, if the multiplication of $X$ is not homotopy associative,
then $r$ is not an $H$-map even though there exists a map between projective planes
extending the suspension of $r$.

Now, the above fact is a special case of the main result of this paper,
which deals with the case that the $H$-space $X$ is an $A_n$-space.
An $A_n$-space is an $H$-space such that the multiplication satisfies higher homotopy associativity
of order $n$.
For example, an $A_2$-space is just an $H$-space, an $A_3$-space is a homotopy associative
$H$-space, and an $A_\infty$-space is a space with the homotopy type of a loop space.

Any $A_n$-space $X$ has an associated space $P_iX$ for each $i$ with $1\le i \le n$ which is
called the projective $i$-space of $X$.
By definition, $P_1X$ is the suspension $\Sigma X$, $P_2X$ is the projective plane,
and $P_\infty X$ is the classifying space of $X$.

Maps preserving $A_n$-space structures are called $A_n$-maps.
An $A_2$-map is an $H$-map, and an $A_\infty$-map is a map homotopic
to a loop map.
See \cite{Iwase-Mimu89a} for the definition.
By definition, if $f \colon X \to Y$ is an $A_n$-map,  then there are maps $P_i f \colon
P_i X \to P_i Y$ $(1\le i \le n)$ such that
\begin{equation}
\label{eqn:qan}
P_1f = \Sigma f, \qquad P_{i+1}f | P_iX \simeq P_if\quad (1\le i \le n-1).
\end{equation}

Then the problem is if the converse of the above fact holds.
To state our result we call a map $f \colon X \to Y$ between $A_n$-spaces a
{\em quasi $A_n$-map} if there are maps $P_i f \colon
P_i X \to P_i Y$ for $(1\le i \le n)$ with \eqref{eqn:qan}.
Then we prove the following

\begin{thm}
\label{theorem}
Let $X$ be an $A_n$-space for some $n\ge 2$.
Then the retraction $r \colon J(X) \to X$ is a quasi $A_n$-map.
\end{thm}

We notice that the above theorem for $n=2$ is just Theorem~\ref{thm:h}.

We can show a fact analogous to Theorem~\ref{thmm} for $A_n$-spaces.
Thus the existence of an $A_{n+1}$-space structure for $X$ is essential 
for the quasi $A_n$-map $r\colon J(X) \to X$ to be an $A_n$-map.
We discuss it in section~\ref{sec:2}.

\section{Proof of the main theorem}

First we recall some facts on the reduced product space given by James \cite{james55a}.
Let $f \colon Z \times J(X) \to Y$ be a map.
Put $f_n=f\circ (id_Z \times \nu_n) \colon Z \times X^n \to Y$ for $n\ge1$,
where $\nu_n \colon X^n \to J(X)$ $(n \ge 1)$ is the canonical map.
Then we have
\[
f_n|Z \times X^{i-1} \times * \times X^{n-i} = f_{n-1}\quad \mbox{for $1\le i \le n$},
\]
where $X^{i-1} \times * \times X^{n-i}$ is identified with $X^{n-1}$ by
the obvious way.

On the other hand, if we have a sequence of maps $(f_n \colon Z \times X^n \to Y)_{n=1,2,\dots}$
with the above property,
then there is a map $f \colon Z \times J(X) \to Y$ such that $f \circ (id_Z \times \nu_n) =f_n$.
Such a sequence $(f_n)_{n=1,2,\dots}$ is called a compatible sequence of invariant maps.

The space $J(X)$ has the homotopy type of $\Omega \Sigma X$.
A homotopy equivalence $s \colon J(X) \to \Omega \Sigma X$ is defined by means of
a compatible sequence of invariant maps $(s_n \colon X^n \to \Omega \Sigma X)_{n=1,2,\dots}$,
where $s_1 \colon X \to \Omega \Sigma X$ is the adjoint of 
$id_{\Sigma X} \colon \Sigma X \to \Sigma X$, and
$s_n$ $(n\ge 2)$ is defined by using the loop multiplication of $\Omega \Sigma X$ as
\[
s_n(x_1,\dots,x_n) = (\cdots (s_1(x_1)\cdot s_1(x_2))\cdots )\cdot s_1(x_n).
\]
Note that to make $(s_n)$ a compatible sequence of
invariant maps we need to modify the loop multiplication so that the constant loop
is the strict unit of the loop multiplication.

Let $e \colon \Sigma \Omega \Sigma X \to \Sigma X$ be the evaluation map,
that is, the adjoint of the $id_{\Omega \Sigma X} \colon \Omega \Sigma X \to \Omega \Sigma X$.
Then we prove the following

\begin{lemma}
\label{lemma}
Let $X$ be an $H$-space and $\varepsilon \colon \Sigma X \to P_2X$ the inclusion.
Then 
$\varepsilon \circ \Sigma r \simeq \varepsilon \circ e \circ \Sigma s$.
\end{lemma}

\begin{proof}
The projective plane $P_2X$ is the mapping cone of the Dold-Lashoff construction
$q \colon X \cup_\mu X \times CX \to \Sigma X$, where
$\mu \colon X \times X \to X$ is the multiplication of $X$.
Morisugi \cite[(1.3)]{morisugi99a} showed that there exists a homotopy
equivalence $X \cup_\mu X \times CX  \to \Sigma (X \wedge X)$ such that if we identify
$X \cup_\mu X \times CX$ with $\Sigma (X \wedge X)$ by this homotopy equivalence,
then $q$ is identified with a map $q' \colon \Sigma (X \wedge X) \to \Sigma X$ with
\[
q' \circ \Sigma \pi\simeq \Sigma p_1+\Sigma p_2- \Sigma \mu \colon \Sigma(X \times X) \to \Sigma X,
\]
where $\pi \colon X \times X \to X \wedge X$ is the quotient map and $p_i$ is the projection to the $i$-th factor.
Thus,
\[
\varepsilon \circ \Sigma \mu \simeq \varepsilon \circ (\Sigma p_1 + \Sigma p_2).
\]

Put $\mu_n=r \circ \nu_n \colon X^n \to X$.
Then $\mu_2=\mu$ and $\mu_n=\mu \circ (\mu_{n-1} \times id)$.
We show that there are homotopies 
$H_n\colon I \times \Sigma X^n \to P_2X$ $(n\ge1)$ between
$\varepsilon \circ \Sigma \mu_n$ and $\varepsilon \circ e \circ \Sigma s_n$ such that
$H_1=\varepsilon \circ p_2$ and
\begin{equation}
\label{eq:H_i}
H_n| I \times \Sigma(X^{j-1} \times * \times X^{n-j}) =H_{n-1}\quad \mbox{for any $1\le j \le n$}.
\end{equation}
Then $(H_n)_{n=1,2,\dots}$ defines a homotopy between $\varepsilon \circ \Sigma r$ and
$\varepsilon \circ e \circ \Sigma s$.

Now $e \circ \Sigma s_2= \Sigma p_1 +\Sigma p_2$ since the adjoint of the both maps are the
same $s_2$.
Thus,
\[
\varepsilon \circ \Sigma \mu_2
\simeq \varepsilon \circ (\Sigma p_1 + \Sigma p_2)
= \varepsilon \circ e \circ \Sigma s_2.
\]
We notice that the above homotopy $H_2 \colon I \times \Sigma X^2 \to P_2X$
can be chosen to be constant of $I \times \Sigma (X \vee X)$.

Let $n>2$.
Suppose inductively that we have $H_i$ for $i<n$ with the desired properties.
Then $H_n$ is defined as the composition of homotopies as follows.
\begin{align*}
\varepsilon \circ \Sigma \mu_n
& = \varepsilon \circ \Sigma \mu \circ \Sigma(\mu_{n-1} \times 1) \\
& \simeq \varepsilon \circ (\Sigma p_1 + \Sigma p_2) \circ \Sigma (\mu_{n-1} \times 1) \\
& = \varepsilon \circ \Sigma \mu_{n-1}\circ \Sigma p' 
 + \varepsilon \circ e \circ \Sigma s_1 \circ \Sigma p_n \\
& \simeq \varepsilon \circ e \circ \Sigma s_{n-1}\circ \Sigma p' 
 + \varepsilon \circ e \circ \Sigma s_1 \circ \Sigma p_n \\
& = \varepsilon \circ e \circ \Sigma s_n,
\end{align*}
where $p' \colon X^n \to X^{n-1}$ is the projection to the first $n-1$-factors,
and the second homotopy is given by using $H_{n-1}$.
It is clear that we can modify $H_n$ to satisfy \eqref{eq:H_i}.
Thus we have $H_n$ for all $n$ by induction.
\end{proof}

Now we prove Theorem~\ref{theorem}.
Theorem~\ref{thm:h} is a special case of Theorem~\ref{theorem}.

\begin{proof}[Proof of Theorem~\ref{theorem}]
Since $J(X)$ is a topological monoid, we have the projective $\infty$-space $P_\infty J(X)$.
It is known that $P_\infty J(X)$ has the homotopy type of $\Sigma X$
such that the inclusion $\Sigma J(X) \to P_\infty J(X)$ followed by
the homotopy equivalence $P_\infty J(X) \simeq \Sigma X$
is homotopic to $e \circ \Sigma s$ (cf. \cite[Proof of Theorem 4.8]{Stas70a}).

Define $P_i r \colon P_iJ(X) \to P_i(X)$ for $2\le i \le n$ by the following composition
\[
P_iJ(X) \subset P_\infty J(X) \simeq \Sigma X \xrightarrow{\varepsilon} P_2X \subset P_iX.
\]
Then by Lemma~\ref{lemma} we have the result.
\end{proof}

\section{$A_n$-form of the retraction}
\label{sec:2}
In this section we show the following theorem which is analogous to Theorem~\ref{thmm}.

\begin{thm}
\label{theorem2}
Let $X$ be an $A_n$-space for some $n\ge 2$.
Then the retraction $r \colon J(X) \to X$ is an $A_{n-1}$-map.
Moreover, if $r$ is an $A_n$-map then the $A_n$-space structure of $X$ is extendable
to an $A_{n+1}$-space structure.
\end{thm}

\begin{proof}
The idea of the proof is not so hard to understand.
But, writing down the explicit proof is very complicated.

Let $\{\mu_i \colon K_i \times X^i \to X\}_{1\le i \le n}$ be the
$A_n$-form on $X$.
The second part of the theorem is a corollary to Iwase-Mimura \cite[P10)]{Iwase-Mimu89a}.
They claim that if $f \colon X \to Y$ and $g \colon Y \to X$ are maps between $A_n$-spaces
such that $g \circ f \simeq id_X$, and
if one of $f$ and $g$ is an $A_n$-map, then the $A_n$-space structure of $X$ is extendable
to an $A_{n+1}$-space structure.
In fact, in our case the extended $A_{n+1}$-form on $X$ is given as follows.

If $\{R_i \colon J_i \times J(X)^i \to X\}_{i\le n}$ is the $A_n$-form on $r$,
then we get $n-1$ higher homotopies
\[
R_n \circ (1 \times \nu_1^s \times \nu_2 \times \nu_1^{n-s-1})\colon
J_n \times X^{n+1} \to X
\quad (1\le s \le n-1).
\]
Then by combining these higher homotopies, we can construct a map
$\mu_{n+1} \colon K_{n+1} \times X^{n+1} \to X$ which extend $\{\mu_i\}_{i\le n}$ to
an $A_{n+1}$-form on $X$.

Next we consider the first part of Theorem~\ref{theorem2}.
An $A_{n-1}$-form $\{R_i \colon J_i \times J(X)^i \to X\}_{2\le i \le n-1}$
is defined by means of compatible sequences of invariant maps
$(R_{i,j} \colon J_i \times J(X)^{i-1}\times X^j \to X)_{j=1,2,\dots}$.

First we define $R_{2,1}$ as the constant homotopy.
For $j\ge2$, $R_{2,j}$ is given as
the composition of $\mu_2 \circ (R_{2,j-1}\times id_X)$
and $\mu_3\circ (1 \times r \times r\circ \nu_{j-1} \times id_X)$.

We don't give the explicit definition for $R_{i,j}$ with $i\ge 3$
since it is almost the same as the case of $R_{2,j}$.
But we just give the following remark.
The homotopy $R_{i,1}$ for $i\ge 3$ can not be a constant homotopy.
For example, $R_{3,1} \colon J_3 \times J(X)^2 \times X \to X$ should be a map
illustrated in Figure~\ref{fig:3},
where the double lines mean constant homotopies.
\begin{figure}
\begin{center}
\caption{$R_{3,1} \colon J_3 \times J(X)^2 \times X \to X$}
\label{fig:3}
\unitlength0.45mm
\footnotesize
\begin{picture}(115,120)(0,-5)
\put(-1,50){\makebox(0,0)[r]{$r(\bmath x)\cdot r(\bmath y \cdot \nu_1(z))$}}
\put(25,101){\makebox(0,0)[rb]{$r(\bmath x \cdot (\bmath y \cdot \nu_1(z)))$}}
\put(90,101){\makebox(0,0)[lb]{$r((\bmath x \cdot \bmath y) \cdot \nu_1(z))$}}
\put(116,50){\makebox(0,0)[l]{$r(\bmath x \cdot \bmath y)\cdot z$}}
\put(25,-1){\makebox(0,0)[rt]{$r(\bmath x) \cdot (r(\bmath y)\cdot z)$}}
\put(90,-1){\makebox(0,0)[lt]{$(r(\bmath x)\cdot r(\bmath y))\cdot z$}}
\put(0,50){\line(1,2){25}}
\put(25,100){\line(1,0){65}}
\put(24,98){\line(1,0){67}}
\put(90,100){\line(1,-2){25}}
\put(87.5,100){\line(1,-2){26}}
\put(0,50){\line(1,-2){25}}
\put(27.5,0){\line(-1,2){26}}
\put(25,0){\line(1,0){65}}
\put(90,0){\line(1,2){25}}
\put(12.5,75){\vector(-1,-2){1}}
\put(14.5,75){\makebox(0,0)[l]{$R_2(t,\bmath x, \bmath y \cdot \nu_1(z))$}}
\put(102.5,25){\vector(-1,-2){1}}
\put(100.5,25){\makebox(0,0)[r]{$R_2(t,\bmath x, \bmath y)\cdot z$}}
\put(57.5,0){\vector(-1,0){1}}
\put(58.5,2){\makebox(0,0)[b]{$\mu_3(t, r(\bmath x), r(\bmath y), z)$}}
\put(57.5,50){\makebox(0,0)[c]{\fbox{$R_{3,1}(\tau, \bmath x, \bmath y, z)$}}}
\end{picture}
\end{center}
\end{figure}
By definition, the homotopy $R_2(t,\bmath x, \bmath y \cdot \nu_1(z))$ is given as the composition
of two homotopies $R_2(t,\bmath x, \bmath y)\cdot z$ and $\mu_3(t, r(\bmath x), r(\bmath y), z)$.
Thus $R_{3,1}$ can be defined by using a suitable degeneracy map 
$\delta \colon J_3 \to J_2$ as
$R_{3,1}(\tau,\bmath x, \bmath y, z)= R_2(\delta(\tau),\bmath x, \bmath y \cdot \nu_1(z))$.
\end{proof}


\begin{thebibliography}{1}

\bibitem{Iwase-Mimu89a}
N.~Iwase and M.~Mimura, \emph{Higher homotopy associativity}, Algebraic
  Topology Proceedings, Arcata 1986 (H.~R.~Miller G.~Carlsson, R. L.~Cohen and
  D.~C. Ravenel, eds.), Lecture Notes in Math., vol. 1370, Springer-Verlag,
  1989, pp.~193--220.

\bibitem{james55a}
I.M. James, \emph{Reduced product spaces}, Ann. of Math. \textbf{62} (1955),
  170--197.

\bibitem{morisugi99a}
K.~Morisugi, \emph{Hopf constructions, {Samelson} products and suspension
  maps}, Contemporary Math. \textbf{239} (1999), 225--238.

\bibitem{Stas61a}
J.~D. Stasheff, \emph{On homotopy abelian ${H}$-spaces}, Cambridge Philos. Soc.
  Proc. \textbf{57} (1961), 734--745.

\bibitem{Stas70a}
\bysame, \emph{${H}$-spaces from a homotopy point of view}, no. 161, Lecture
  Notes in Math., Springer-Verlag, 1970.

\end{thebibliography}
\end{document}